\documentclass{amsart}
\usepackage{amssymb, mathtools, fullpage, color,cite,booktabs}
\usepackage[colorinlistoftodos]{todonotes}
\usepackage{siunitx}
\usepackage{url}

\theoremstyle{definition}

\newenvironment{red}{\relax\color{red}}{\relax}
\newenvironment{blue}{\relax\color{blue}}{\hspace*{.5ex}\relax}

\newcommand{\ber}{\begin{red}}
\newcommand{\er}{\end{red}}
\newcommand{\beb}{\begin{blue}}
\newcommand{\eb}{\end{blue}}

\newcommand{\botrule}{\midrule}

\numberwithin{equation}{section}

\begin{document}

\title[Murmurations: a case study in AI-assisted mathematics]{Murmurations: a case study in AI-assisted mathematics}

\date{\today}

\author[Y.-H. He]{Yang-Hui He$^{\dagger}$}
\thanks{$^{\dagger}$Supported in part by UK STFC grant ST/J00037X/3.}
\address{London Institute for Mathematical Sciences, Royal Institution, London W1S 4BS, UK \hfill \break \indent   Merton College, Oxford, OX14JD, UK}
\email{yh@lims.ac.uk}

\author[K.-H. Lee]{Kyu-Hwan Lee}
\address{Department of Mathematics, University of Connecticut, Storrs, CT       06269, USA  \hfill \break \indent Korea Institute for Advanced Study, Seoul     02455, Republic of Korea}
\email{khlee@math.uconn.edu}

\author[T. Oliver]{Thomas Oliver}
\address{University of Westminster, 115 New Cavendish Street, London W1W 6UW, UK}
\email{T.Oliver@westminster.ac.uk}

\author[A. Pozdnyakov]{Alexey Pozdnyakov}
\address{Department of Mathematics, Fine Hall A02, Washington Road,
Princeton University, NJ 08544-1000 USA}
\email{ap5763@princeton.edu}

\begin{abstract}
We report the emergence of a striking new phenomenon in arithmetic, which we call \emph{murmurations}. 
First observed experimentally through averages over large arithmetic datasets, murmurations can be detected and analyzed using standard interpretability tools from machine learning, including principal component weightings, saliency curves, and convolutional filters.
Although discovered computationally, they constitute a genuinely new and intriguing phenomenon in arithmetic that can be formulated and investigated using established tools of number theory. 
In particular, murmurations encode subtle information about Frobenius traces and naturally belong to the framework of arithmetic statistics. 
More precisely, murmurations connect to central themes surrounding the conjecture of Birch and Swinnerton-Dyer and perspectives from random matrix theory. In this paper, we present an overview of murmurations, contextualizing them within number theory and AI.
\end{abstract}

\maketitle

\newcommand{\comment}[1]{}

\section{Introduction}\label{sec:introduction}

Prime numbers and elliptic curves are among the central objects of number theory. 
Although both have been studied for thousands of years, they continue to reveal surprising patterns. 
Modern number theory often investigates these objects statistically: instead of studying a single prime or a single curve, one examines large datasets and searches for patterns in aggregate behaviour. 
A famous example is the conjecture of Birch and Swinnerton-Dyer (BSD), which grew out of computer experiments dating back to the 1950s and relates statistical/analytic properties of elliptic curves and their associated $L$-functions to deep algebraic invariants of the curves. 
At the heart of these questions lie Frobenius traces, which are integers that measure how an elliptic curve behaves modulo a prime. 
Understanding the statistical distribution of these traces is a major theme in arithmetic statistics and connects with ideas from random matrix theory.

In this paper, we describe the discovery of a new and unexpected pattern for Frobenius traces, which we call {\em murmurations} \cite{He:2022pqn}. 
What makes this discovery particularly striking is that elliptic curve data is already extremely well studied: large databases such as the LMFDB have been carefully compiled and analyzed by experts for many years \cite{lmfdb}. 
Nevertheless, by applying relatively simple techniques, we observed a striking, scale-invariant, oscillatory pattern that had previously gone unnoticed. 
More sophisticated machine learning tools, including saliency curves and convolutional filters, also detect murmurations \cite{pozdnyakov2024predicting, big-group}, 
however, their signals are ultimately dominated by the well-known Mestre–Nagao sums, which are classical quantities that are already understood to correlate strongly with rank predictions. 
This contrast highlights both the subtlety of the new phenomenon and the importance of carefully interpreting machine learning outputs.

\medskip

The study of prime numbers has long had a statistical flavour. 
Although Euclid proved over two thousand years ago that there are infinitely many primes, mathematicians soon began asking quantitative questions, such as how many primes lie below a given bound. 
In the late eighteenth century, Gauss and Legendre, guided by extensive numerical tables, independently proposed an approximation for counting primes. 
Their conjecture was eventually proved in 1896 using tools from complex analysis, following crucial insights of Riemann in 1859. 
Riemann’s hypothesis about the zeros of his zeta function, still unproven today, remains one of the central problems in mathematics \cite{carlson2006millennium}.
Like the conjecture of Birch and Swinnerton-Dyer, the Riemann hypothesis illustrates a recurring theme: large-scale numerical patterns in prime-related data often point toward profound underlying structure.

The murmurations phenomenon therefore sits within a broader story about statistics and the arithmetic information they encode. 
While the pattern was first detected through computational and data-driven methods, it can be formulated in purely mathematical terms and studied using established tools of analytic and algebraic number theory. 
Our goal in this article is to present an accessible overview of murmurations from both arithmetic and AI perspectives, explaining how a simple experimental approach led to a genuinely new question in a classical and intensely studied area of mathematics.
The discovery of murmurations forms part of a broader and rapidly developing interaction between artificial intelligence and pure mathematics. 
In recent years, machine learning has been used to assist in conjecture formation, pattern recognition, and large-scale data analysis across several areas of mathematics and theoretical physics \cite{He:2017aed,He:2018jtw,Krefl:2017yox,Carifio:2017bov,Ruehle:2017mzq,lample2019deep,udrescu2020ai,iten2020discovering,He:2021oav,davies2021advancing,raayoni2021generating,krenn2022scientific,fawzi2022discovering,ML2023,weinan2021dawning,he2019learning,he2020graph, CKVcomms, CKVneura, Lee2024MLKronecker}. 
One influential line of work proposes a systematic framework in which a mathematical hypothesis is formulated, data are generated, a supervised learning model is trained, and attribution methods are then used to suggest new conjectures~\cite{davies2021advancing}. 
This approach has demonstrated that machine learning can successfully identify meaningful mathematical structure~\cite{He:2021oav}.

Our discovery differs in both method and outcome to the framework presented in~\cite{davies2021advancing}. 
Murmurations emerged from exploratory data analysis applied to arithmetic statistics. 
The key step was not the training of a highly complex model, but the careful choice of how to represent number-theoretic data and how to average it. The resulting oscillatory behaviour was not an obvious refinement of known results, nor something that could be quickly verified or dismissed. On the contrary, its meaning is still being investigated, and it has already prompted new theorems and ongoing theoretical work \cite{cowan2023murmurations, ICERM, Z23, Zu, BBLL, LOPdir, LDMTF,BBLLDSHZ, IERCMNS, cowan2024murmurations, wang, mppower, MHRZI, MHRZII, KCLDL}. In this sense, murmurations illustrate a less linear and less predictable interaction between AI and mathematics than is sometimes portrayed.
A comprehensive exploration of the scope of the discovery will appear in joint work of the present authors with Sutherland \cite{HLOPS}.
The non-triviality contrasts with \cite{udrescu2020ai,davies2021advancing,fawzi2022discovering}, in which the 
results proposed by AI could be checked or proven in a relatively short time, and also with \cite{He:2017aed,he2019learning} where the patterns spotted by AI are difficult to turn into a precise mathematical conjecture.

The discovery of murmruations belongs to a developing programme of AI-driven number theory \cite{He:2020kzg,He:2020tkg,He:2022pqn,He:2020qlg,Amir:2022sab,HLOPS}, which builds on some pioneering investigations in the same domain \cite{raayoni2021generating, He:2018jtw, Alessandretti:2019jbs}. In this setting, arithmetic objects (such as elliptic curves) are encoded as high-dimensional data via the coefficients of their associated $L$-functions. 
These coefficients include the Frobenius traces discussed earlier, which play a central role in conjectures such as Birch and Swinnerton–Dyer. 
Large curated databases, most notably the LMFDB, provide an unprecedented supply of structured arithmetic data. 
Machine learning methods are particularly well suited to detecting patterns in such high-dimensional environments, where purely human inspection would be infeasible.

At the same time, experience shows that blindly applying sophisticated algorithms to large databases is rarely effective. 
For example, advanced interpretability tools such as saliency curves and convolutional filters can indeed detect signals in elliptic curve data, but their outputs are often dominated by well-known quantities such as the Mestre–Nagao sums. 
Extracting genuinely new information requires mathematical insight into how the data are organised and which features are likely to encode deep arithmetic behaviour. 
The discovery of murmurations therefore did not arise from automation alone, but from a hybrid approach: human expertise guided the representation of the data and the design of the experiment, while machine learning tools amplified our ability to detect subtle aggregate structure.

This interplay reflects a broader lesson. 
While AI excels at large-scale pattern recognition, the formulation of meaningful mathematical questions (and the interpretation of the answers) remains deeply dependent on human understanding \cite{buzzard,He:2021oav,williamson2023deep,He:2018jtw,davies2021advancing,mishra2023mathematical,kolpakov2023impossibility}. 
Murmurations provide a case study in how classical number theory, arithmetic statistics, and interpretable machine learning can interact to produce genuinely new mathematics.

\section{Murmurations: an arithmetic perspective}
Life is not fair.
Perhaps the most dramatic example of this old adage is that even arithmetic is \textsl{biased}.
In this section, we will offer examples within the context of prime numbers and elliptic curves.

\subsection{Biases for prime numbers}
Consider a prime number $p$ other than 2, the smallest and only even prime number.
If we divide $p$ by 4, we will get a remainder $1$ or $3$ because $p$ is odd.
In 1853, Chebyshev \cite{chebyshev1853lettre} noticed that, if one fixes an upper bound and considers all the odd primes less than it, then, when dividing by $4$, a remainder of $3$ seems to occur {\it more often} than a remainder of $1$.

To illustrate this bias, in Table~\ref{tab1}, we list all odd primes less than 100, and record the remainder when divided by $4$.
The average of the possible remainders is $\frac12(1+3)=2$, and we record the discrepancy from this average for each prime.
Thus a remainder of 1 constitutes a positive discrepancy, and a remainder of 3 constitutes a negative discrepancy.
In addition, we record the aggregate discrepancy, which is the sum of the discrepancies for all of the primes until that point. 
We notice that the aggregate discrepancy is \textsl{negatively biased}.
In fact, in our table, the aggregate discrepancy does not ever take a positive value. 
In other words, there are always at least as many occurrences of the remainder $3$ as there are of the remainder $1$.

\begin{center}
\begin{table}[h]
\begin{center}
\caption{Remainders of odd primes when divided by $4$. }\label{tab1} 
{\tiny \begin{tabular}{@{}rrrr|rrrr@{}}
\toprule
$\begin{array}{c}\text{Odd}\\ \text{prime}\end{array}$ & Remainder & Discrepancy& $\begin{array}{c}\text{Aggregate}\\ \text{discrepancy}\end{array}$ & $\begin{array}{c}\text{Odd}\\ \text{prime}\end{array}$ & Remainder & Discrepancy & $\begin{array}{c}\text{Aggregate}\\ \text{discrepancy}\end{array}$  \\
\midrule
$3   $&$ 3 $&$ -1$& $-1$  &
$ 43 $&$ 3 $&$-1 $& $-1$ \\$5    $&$ 1 $&$1  $& $0$ &
$ 47 $&$ 3$&$-1 $ & $-2$\\$7    $&$ 3 $&$-1 $& $-1$ &
$ 53 $&$ 1 $&$1 $  & $-1$\\$11   $&$ 3 $&$-1 $& $-2$ &
$ 59 $&$ 3$&$-1 $ & $-2$\\$13    $&$ 1$&$1 $& $-1$ &
$ 61 $&$ 1$&$1   $ & $-1$\\$17    $&$ 1 $&$1 $& $0$ &
$ 67 $&$ 3 $&$-1$   & $-2$\\$19   $&$ 3 $&$-1 $& $-1$ &
$ 71 $&$ 3$&$-1 $ & $-3$\\$23    $&$ 3$&$-1 $& $-2$ &
$ 73 $&$ 1 $&$1 $  & $-2$\\$29    $&$ 1 $&$1 $& $-1$ &
$ 79 $&$ 3 $&$-1$  & $-3$\\$31   $&$ 3 $&$-1 $& $-2$ &
$ 83 $&$ 3$&$-1 $& $-4$\\$37    $&$ 1$&$1 $& $-1$ &
$ 89 $&$ 1 $&$1  $ & $-3$\\$41    $&$ 1 $&$1 $& $0$ & 
$ 97 $&$ 1$&$1  $ & $-2$\\
\botrule 
\end{tabular} }
\end{center}
\end{table}
\end{center}

\begin{figure}[h] 
\centering
\includegraphics[width=0.9\textwidth]{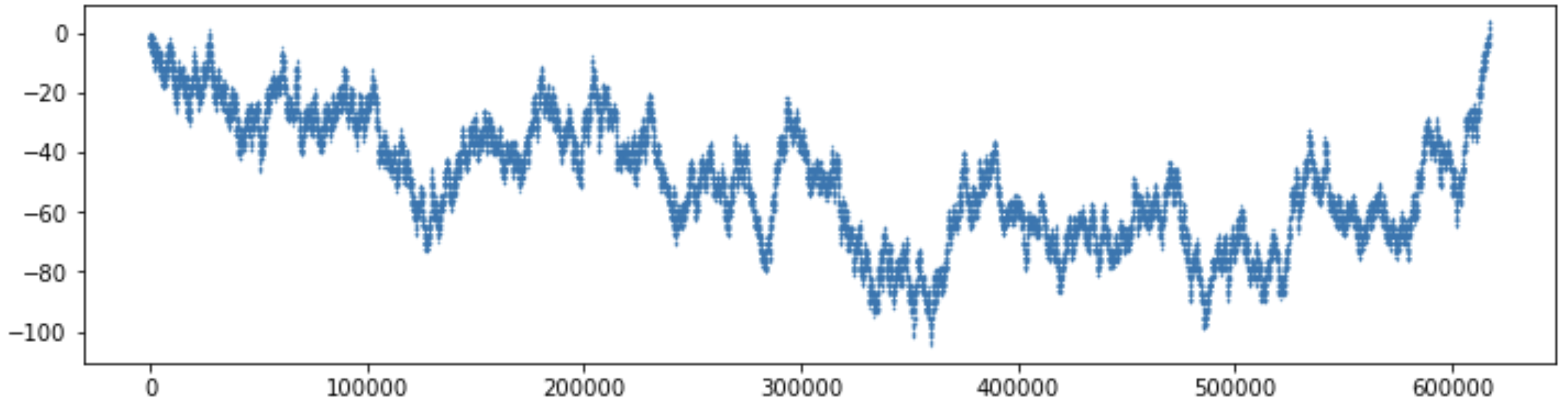}
\caption{A plot of the primes ($x$-axis) against the aggregate discrepancy ($y$-axis).
}\label{fig1}
\end{figure}

One can ask if, were we to count high enough, could the count of 1 ever exceed the count of 3, that is, is the aggregate discrepancy ever positive?
Perhaps surprisingly, the answer is yes, though we need to wait a while to see it. 
The first prime for which the aggregate discrepancy is positive is $\num{26861}$, and, after that, the next is $\num{616841}$.
Whilst it may seem extremely rare, there are infinitely many primes for which the aggregate discrepancy is positive and, under the additional assumption of the Riemann Hypothesis, we can, in some sense, quantify \textsl{exactly} how often this occurs \cite{rubinstein1994chebyshev}.
As a consequence, we can reasonably assert that the negative bias observed in Table~\ref{tab1} persists for \textsl{most} $p$.
This observation is represented graphically in Figure~\ref{fig1}.

\subsection{Elliptic Curves}
\begin{figure}[h] 
\centering
\includegraphics[scale=0.3]{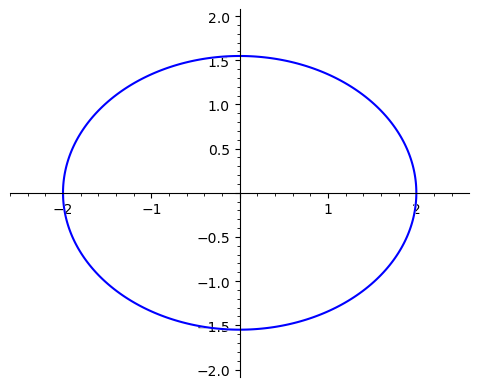} \hskip 1 cm
\includegraphics[scale=0.3]{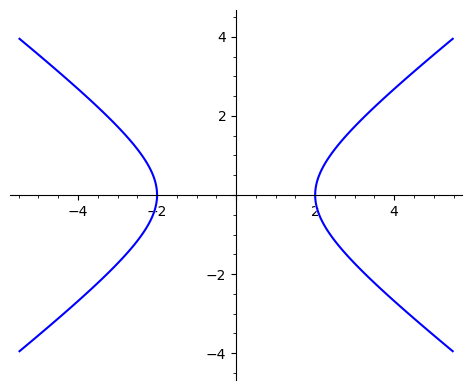}
\caption{Graphs for conic sections defined by the quadratic equations (left) $3x^2+5y^2=12$ and (right) $3x^2-5y^2=12$.}\label{fig-conic}
\end{figure}

\begin{figure}[h] 
\centering
\includegraphics[scale=0.3]{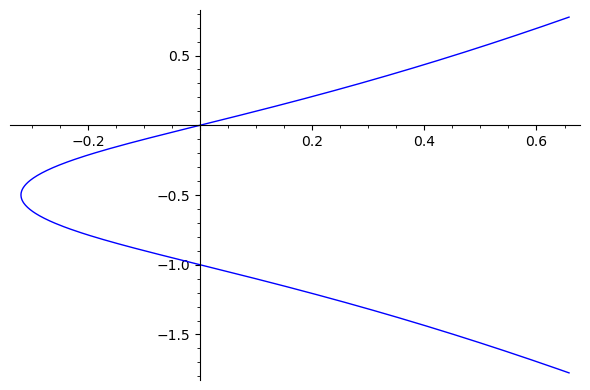} \hskip 1 cm
\includegraphics[scale=0.3]{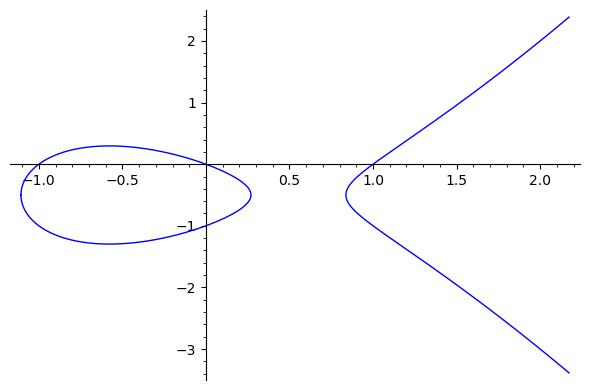}
\caption{Graphs for elliptic curves defined by the cubic equations (left) $y^2+y=x^3+x^2+x$  and (right)  $y^2+y=x^3-x$.}\label{fig-elliptic-curve}
\end{figure}

Elliptic curves are central to number theory and feature in everything from the aforementioned BSD conjecture \cite{carlson2006millennium} to Wiles' proof of Fermat's Last Theorem \cite{wiles1995modular}.
Elliptic curves are planar curves described by a special cubic equation in two variables $x$ and $y$ (see examples in Figure \ref{fig-elliptic-curve}).
To the uninitiated, it may seem that elliptic curves are only marginally more complicated than the conic sections defined by \textsl{quadratic} equations in two variables (see examples in Figure \ref{fig-conic}).
Nevertheless, the cubic variable ensures that the situation is much more challenging.
So challenging, in fact, that these equations underpin crucial aspects of modern cybersecurity \cite{nistECC}.

The plots in Figures \ref{fig-conic} and \ref{fig-elliptic-curve} are the \textsl{real} loci studied by Apollonius (200BC) and Newton \cite{harris1725lexicon}.
In number theory, instead of considering solutions $(x,y)$ with real co-ordinates as in classical geometry, the focus is on $(x,y)$ being integers and rational numbers (ratios of integers).
Finding such rational solutions to polynomial equations is notoriously difficult and is known as the {\it Diophantine problem}, first considered systematically by Diophantus (200AD).

When seeking to find rational solutions to a polynomial in two variables there are three possibilities, namely, there could exist (1) no solutions, (2) finitely many solutions, or (3) infinitely many solutions.
In the case of conic sections, it is known that, except for the degenerate case of a single solution, only cases (1) and (3) can occur.
One reason that elliptic curves are interesting is that they are the simplest context in which case (2) of finitely many, but more than 1, is realised.
On the other hand, for more complicated curves\footnote{The complexity we are referring to here is known as the \textsl{genus}.} it is known that the possibility (3) never occurs.
In short, elliptic curves present an interesting threshold for the Diophantine problem.

Let us consider some examples.
The elliptic curve defined by $y^2+y=x^3+x^2+x$ (see Figure~\ref{fig-elliptic-curve}) has only finitely many rational solutions, namely, the two integer solutions $(0, 0)$ and $(0, -1)$.
On the other hand, the elliptic curve defined by $y^2+y=x^3-x$ (see Figure~\ref{fig-elliptic-curve}) has 10 integer solutions:
\[(-1, 0), \ \ (-1, -1), \ \ (0, 0), \ \ (0, -1), \ \ (1, 0), \ \ (1, -1), \ \ (2, 2), \ \ (2, -3), \ \ (6, 14), \ \  (6, -15),\] 
but infinitely many other rational solutions, some of which we list below: 
\[\left(\tfrac 1 4 , -\tfrac 5 8 \right), \ \  \left(-\tfrac 5 9  , \tfrac 8 {27} \right), \ \ \left(\tfrac{21}{25}  , -\tfrac{69}{125}  \right), \ \  \left(-\tfrac{20}{49} , -\tfrac {435}{343} \right), \ \ \left(\tfrac{161}{16} , -\tfrac{2065}{64} \right), \ \  \left(\tfrac{116}{529}, -\tfrac{3612}{12167} \right), \dots\] 
If an elliptic curve has only finitely many rational solutions, then we say it has \textsl{rank} 0.
For example, the left curve in Figure~\ref{fig-elliptic-curve} has rank $0$.
For curves with infinitely many rational solutions, we assign a positive integer rank reflecting the ``size'' of the set of rational points.
For example, the right curve in Figure~\ref{fig-elliptic-curve} has rank $1$.

The terminology \textsl{rank} used here is taken from group theory in abstract algebra.
The famous BSD conjecture relates this algebraic concept of rank to an entirely different formulation that uses complex analysis.
It is this link between the seemingly disjoint worlds of analysis and algebra that gives the BSD conjecture its depth.

We will not give a precise definition of the rank here; the reader may refer to \cite[Section~1]{RS} for details. 
Although the rank is rigorously defined, there is no known algorithm to compute it for an arbitrary elliptic curve \cite[Section~3]{RS}. 
Moreover, it is not known which positive integers can occur as the rank of an elliptic curve.
A related invariant is the \textsl{root number}. 
It equals $1$ for elliptic curves of rank $0$ and $-1$ for those of rank $1$. 
More generally, it is conjectured that the root number is $1$ for curves of even rank and $-1$ for curves of odd rank. 
This statement is known as the \emph{parity conjecture} \cite{Dok}.

\subsection{Biases for elliptic curves}

A standard trick for solving a polynomial equation with integers is to consider its reduction modulo a prime number $p$.
What this means is that, rather than being equal, both sides of the equation have the same remainder when divided by $p$.
For example, the elliptic curve defined by $y^2+y=x^3-x$ has a solution modulo $p=11$ given by $(x,y)=(7,2)$. 
This is because $x^3-x=7^3-7=336$ and $y^2+y=2^2+2=6$ both have remainder $6$ when divided by $11$. 
Since the remainder when dividing by $11$ must be an integer between $0$ and $10$, there can be at most $11\times11=121$ solutions to a two-variable equation modulo 11.
By exhaustively checking all possibilities, one can see that  $y^2+y=x^3-x$ has 16 solutions modulo $11$. 

\begin{figure}[h] 
\centering
\includegraphics[scale=0.4]{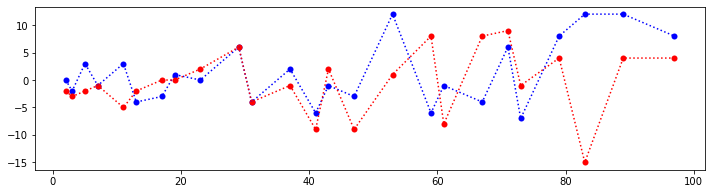}
\includegraphics[scale=0.4]{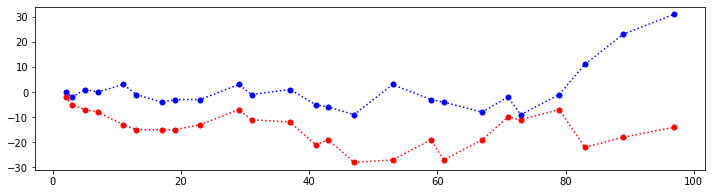}
\caption{(Top) A plot showing how the discrepancy varies, with $y^2+y=x^3-x$ in red and $y^2+y=x^3+x^2+x$ in blue. (Bottom) A plot showing how the aggregate discrepancy varies for the same curves.} \label{fig-ec-dev}
\end{figure}

Now let us work abstractly and consider a generic elliptic curve $E$ modulo a generic prime number $p$. 
We will denote the number of solutions to $E$ modulo $p$ (plus an additional solution ``at infinity'') by $\#E_p$.
There is an efficient algorithm to compute this number exactly.
Just like in the previous section, in the equation below we introduce a sort of discrepancy for the solution counts\footnote{In geometric terms, the discrepancy quantifies the difference between $E$ and a straight projective line. This accounts for the ``$p+1$'' in the displayed equation.}.
It is customary to use the following notation  for this discrepancy:
\[a_p(E) := p + 1 - \#E_p.\]
The so-called \textsl{Hasse bound} implies that the discrepancy does not grow too large as $p$ increases.
In Tables \ref{tab-single-elliptic-1} and \ref{tab-single-elliptic}, we list $a_p(E)$ for our two example curves, for all the primes up to 100 (the same information is presented in graphical form in Figure \ref{fig-ec-dev}). 
In Table \ref{tab-single-elliptic} we see a negative bias in the aggregate discrepancy, while in Table \ref{tab-single-elliptic-1} the sign changes several times. 
The presence of negative aggregates in Table \ref{tab-single-elliptic} reflects the widely accepted heuristic that larger rank curves should have more solutions modulo primes \cite{Martin-Pharis}.
We will revisit this belief in the next section.

An ingenious idea is to gather together the discrepancies for each prime into a special function, known as the (elliptic) $L$-function.
When this is done, the $a_p(E)$ appear as  \textsl{Dirichlet coefficients}.
It is the elliptic $L$-function that is central both to the proof of Fermat’s Last Theorem and to the formulation of the Birch–Swinnerton-Dyer conjecture. Furthermore, 
$L$-functions are a cornerstone of the Langlands programme, which has become a central theme of much modern mathematical research.

\begin{center}
\begin{table}[h]
\begin{center}
\caption{Counting solutions on the fixed elliptic curve $y^2+y=x^3+x^2+x$, which has rank 0.}\label{tab-single-elliptic-1}
\begin{tabular}{@{}rrrr|rrrr@{}}
\toprule
Prime& Count& Discrepancy & $\begin{array}{c}\text{Aggregate}\\ \text{discrepancy}\end{array}$ & Prime& Count& Discrepancy & $\begin{array}{c}\text{Aggregate}\\ \text{discrepancy}\end{array}$   \\
\midrule
$2$ & $2$ & $0$ & $0$  &
$43$ & $44$ & $-1$ & $-6$ \\
$3$ & $5$ & $-2$ & $-2$  &
$47$ & $50$ & $-3$ & $-9$\\
$5$ & $2$ & $3$ & $1$ &
$53$ & $41$ & $12$ & $3$\\
$7$ & $8$ & $-1$ & $0$ &
$59$ & $65$ & $-6$ & $-3$\\
$11$ & $8$ & $3$ & $3$ &
$61$ & $62$ & $-1$ & $-4$\\
$13$ & $17$ & $-4$ & $-1$ &
$67$ & $71$ & $-4$ & $-8$\\
$17$ & $20$ & $-3$ & $-4$ &
$71$ & $65$ & $6$ & $-2$\\
$19$ & $18$ & $1$ & $-3$ &
$73$ & $80$ & $-7$ & $-9$\\
$23$ & $23$ & $0$ & $-3$ &
$79$ & $71$ & $8$ & $-1$\\
$29$ & $23$ & $6$ & $3$ &
$83$ & $71$ & $12$ & $11$\\
$31$ & $35$ & $-4$ & $-1$ &
$89$ & $77$ & $12$ & $23$\\
$37$ & $35$ & $2$ & $1$ &
$97$ & $91$ & $8$ & $31$\\
$41$ & $47$ & $-6$ & $-5$ &&&&\\ 
\botrule
\end{tabular}
\end{center}
\end{table}
\end{center}

\begin{center}
\begin{table}[h]
\begin{center}
\caption{Counting solutions on the fixed elliptic curve $y^2+y=x^3-x$, which has rank 1. }\label{tab-single-elliptic}
\begin{tabular}{@{}rrrr|rrrr@{}}
\toprule
Prime& Count& Discrepancy & $\begin{array}{c}\text{Aggregate}\\ \text{discrepancy}\end{array}$ & Prime& Count& Discrepancy & $\begin{array}{c}\text{Aggregate}\\ \text{discrepancy}\end{array}$  \\
\midrule
$2$ & $4$ & $-2$ & $-2$  &
$43$ & $41$ & $2$ & $-19$ \\
$3$ & $6$ & $-3$ & $-5$  &
$47$ & $56$ & $-9$ & $-28$\\
$5$ & $7$ & $-2$ & $-7$ &
$53$ & $52$ & $1$ & $-27$\\
$7$ & $8$ & $-1$ & $-8$ &
$59$ & $51$ & $8$ & $-19$\\
$11$ & $16$ & $-5$ & $-13$ &

$61$ & $69$ & $-8$ & $-27$\\
$13$ & $15$ & $-2$ & $-15$ &
$67$ & $59$ & $8$ & $-19$\\
$17$ & $17$ & $0$ & $-15$ &
$71$ & $62$ & $9$ & $-10$\\
$19$ & $19$ & $0$ & $-15$ &
$73$ & $74$ & $-1$ & $-11$\\
$23$ & $21$ & $2$ & $-13$ &
$79$ & $75$ & $4$ & $-7$\\
$29$ & $23$ & $6$ & $-7$ &
$83$ & $98$ & $-15$ & $-22$\\
$31$ & $35$ & $-4$ & $-11$ &
$89$ & $85$ & $4$ & $-18$\\
$37$ & $38$ & $-1$ & $-12$ &
$97$ & $93$ & $4$ & $-14$\\
$41$ & $50$ & $-9$ & $-21$ &&&&\\ 
\botrule
\end{tabular}
\end{center}
\end{table}
\end{center}

\subsection{Datasets of elliptic curves}

The breakthrough in \cite{He:2022pqn} involves a different sort of aggregate behaviour, which somehow defies the expectation that higher rank curves tend to have more solutions mod $p$.
Rather than sum over prime numbers, we instead sum over \textsl{similar} curves.
Let us first explain what this means.

One relevant notion is that of \textsl{isogeny}.
Two elliptic curves are said to be \textsl{isogenous} if they have the same $a_p$ for every $p$. 
This can happen for rather different looking equations, for example, $y^2+y=x^3-x$ is isogenous to $y^2+y=x^3-x^2-7820x-263580$ and $y^2+y=x^3-x^2-10x-20$.
Since our primary interest is the Frobenius traces $a_p$, we will treat a collection of isogenous curves as a single object, known as an \textsl{isogeny class}.
For example, the isogeny class for $y^2+y=x^3-x$ contains the three curves listed previously.
Isogenous curves have many properties in common, for example, they have the same rank.

\begin{figure}[h] 
\centering
\includegraphics[scale=0.4]{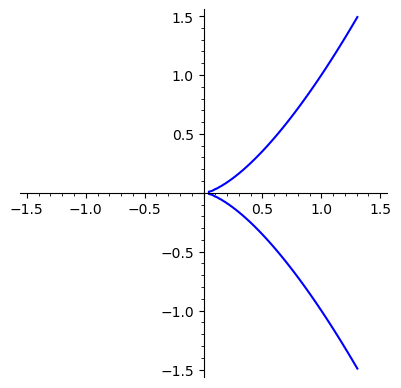} \hskip 1 cm
\includegraphics[scale=0.4]{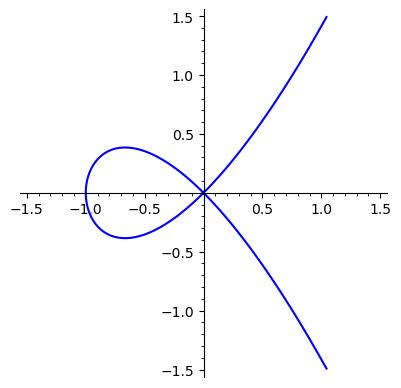}
\caption{(Left) The curve $y^2=x^3$ has a cuspidal singularity, and (right) the curve $y^2=x^3+x^2$ has a nodal singularity.} \label{fig-cusp-node}
\end{figure}

We will order the isogeny classes by \textsl{conductor}, which is a measure of complexity.
Reducing a polynomial equation can fundamentally change the geometry by introducing \textsl{singularities}.
The conductor $N(E)$ of an elliptic curve $E$ keeps track of this: it is an integer divisible only by the primes producing a singular reduction.
The reductions of elliptic curves may exhibit two types of singularity, namely, cusps and nodes. 
These singularities are depicted in Figure \ref{fig-cusp-node}.
The type of singularity induced by a prime determines the power of that prime in the factorization of the conductor.
It is also possible to take an analytic perspective on the conductor, which measures the density of zeros for the corresponding $L$-function.

Bringing together various concepts introduced above, we can explain what we mean by similar curves. 
Namely, isogeny classes with the same rank parity and conductors within a given range of values.
Given an interval $I$ on the real number line and a rank parity $\varepsilon \in \{ 0,1 \}$, we may write collections of similar curves symbolically as a \textsl{dataset}:
\[S(\varepsilon,I) := \left\{
\text{isogeny classes }E :
N(E) \in I,
\ \mathrm{rank}(E)\equiv \varepsilon \bmod2
\right\},\]
where we have explicitly labeled $\varepsilon$ and $I$ because this set depends on both choices.

\subsection{Murmurations of elliptic curves}

Our key discovery of the murmuration phenomenon refers to an aggregate behaviour for similar elliptic curves.
More precisely, for an interval $I$, we average the values $a_p(E)$ over all isogeny classes $E$ in $S(\varepsilon,I)$.
This yields a function on the set of primes.
Mathematically, we write
\[
m_{\varepsilon, I}(p) :=
\frac{1}{\vert S(\varepsilon,I) \vert} \sum\limits_{E \in S(\varepsilon,I)} a_p(E).
\]
In Table \ref{tab4} we provide some numerical values for this function for a given interval.
In Figure \ref{fig-mur}, we plot the analogous graphs for larger intervals.
We contrast the red points in Figure~\ref{fig-mur} with those in Figure~\ref{fig-ec-dev} (bottom image).
In the latter (which is concerned with a single rank 1 curve) the red points are consistently below the $x$-axis (as one expects).
In the former (which is concerned with families of similar odd-rank curves) the red points oscillate around the $x$-axis (which is a surprise).

Whilst the oscillation alone is unexpected, what is more surprising is the \textsl{scale invariance}.
This is explored in more detail in \cite{HLOPS}.
In each frame of Figure~\ref{fig-mur}, we have plotted the average value of the discrepancy for primes up to the start of the interval.
We have then scaled all images uniformly so that the $x$-axes have the same length.
We not only observe that the oscillation persists for different intervals, but that the key features (peaks and crossover points) occur at the same locations. 
This is despite the fact that the datasets involved in the average are completely disjoint!

\vspace{5pt}

\begin{table}[h]
\begin{center}
\caption{Numerical value for the average of $a_p(E)$, as $E$ varies over isogeny classes with conductor in $[200,400]$ for each rank parity. There are 256 isogeny classes with even rank, and 190 with odd rank. }\label{tab4} 
\begin{tabular}{@{}rrr|rrr@{}}
\toprule
Prime & $\begin{array}{c} \text{Average over}\\ \text{even ranks}\end{array}$ & $\begin{array}{c} \text{Average over}\\ \text{odd ranks}\end{array}$  & Prime & $\begin{array}{c} \text{Average over}\\ \text{even ranks}\end{array}$  & $\begin{array}{c} \text{Average over}\\ \text{odd ranks}\end{array}$  \\
\midrule
2   &$0.137$ &$-0.242$ & 43 &$-0.953$  &$-0.321$  \\
3    &$0.266$  &$-0.605$ & 47 &$-1.148$ &$0.332$  \\
5    &$0.715$ &$-1.100$ & 53 &$-1.484$  &$1.505$  \\
7   &$0.871$ &$-1.579$ & 59 &$-1.062$  &$0.668$ \\
11    &$1.078$  &$-1.937$ & 61 &$-1.535$  &$0.968$   \\
13    &$0.770$ &$-2.263$  & 67 &$-1.648$  &$0.316$   \\
17   &$0.906$  &$-1.905$ & 71 &$-0.523$  &$0.463$  \\
19    &$0.777$ &$-1.963$  & 73 &$-1.039$  &$-0.342$  \\
23    &$0.918$ &$-1.321$  & 79 &$-0.742$  &$-0.542$   \\
29   &$0.039$  &$-0.974$ & 83 &$0.031$  &$-0.189$ \\
31    &$-0.266$ &$-1.047$ & 89 &$0.410$  &$-1.379$  \\
37    &$-0.863$  &$-0.684$ & 97 &$0.703$  &$-1.789$  \\
41 &$-0.812$&$-0.105$&&&\\
\botrule
\end{tabular}
\end{center}
\end{table}

\comment{
\begin{figure}[h] 
\centering
\includegraphics[scale=0.4]{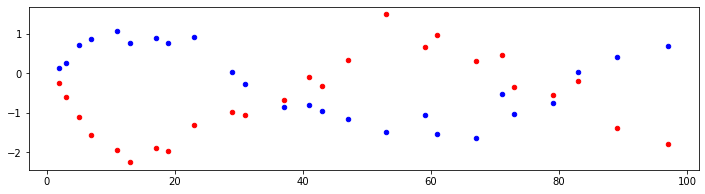}
\caption{Plot of the average of $a_p(E)$, as $E$ varies over isogeny classes with conductor in $[200,400]$ for each rank parity. The even rank curves are in blue, and the odd rank curves are in red.} \label{fig-mur-small}
\end{figure}
}

\vspace{5pt}

\begin{figure}[h] 
\centering
\includegraphics[scale=0.4]{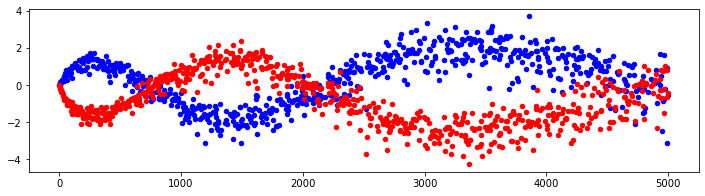}
\includegraphics[scale=0.4]{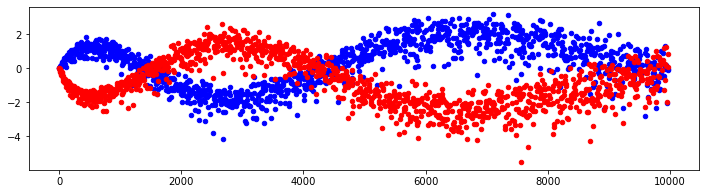}
\includegraphics[scale=0.4]{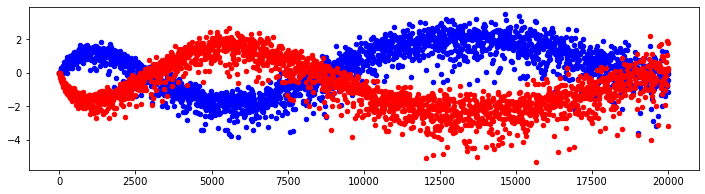}
\caption{
Plot of the average of $a_p(E)$, as $E$ varies over isogeny classes with a fixed rank parity and conductor in (Top) $[5000,10000]$, (Middle) 
$[10000,20000]$, (Bottom)
$[20000,40000]$. In each image, the even rank curves are in blue, and the odd rank curves are in red.
} \label{fig-mur}
\end{figure}

\vspace{5pt}

\section{Murmurations: an AI perspective}
In the introduction, we referred to a programme on \textsl{AI-driven number theory} which was launched in \cite{He:2020kzg,He:2020tkg,He:2022pqn,He:2020qlg,Amir:2022sab,HLOPS}. Based on the exposition so far, which is expressed in purely number theoretic terms, the reader might reasonably question what role, if any, was played by AI.
In fact, the presented formulation of murmurations as a straightforward average stems from a remarkable success in interpreting the outputs of machine learning algorithms.
In what follows, it is our intention to elucidate the way in which machine learning pointed us towards murmuration.
More precisely, we will overview evidence suggested by two paradigms of machine learning (namely, unsupervised and supervised), both in general terms, and specifically for our data.

\begin{figure}[h]
\centering
\includegraphics[scale=0.3]{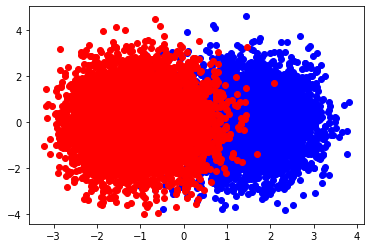} \hskip 1 cm
\includegraphics[scale=0.3]{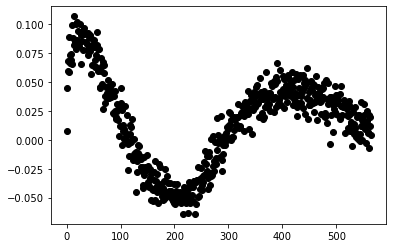}
\caption{PCA on point cloud of elliptic curves with conductor in $[5000, 10000]$ (represented as 564-dimensional vectors; this is chosen because there are 564 primes $<2^{12}$) with resulting projection on the left and vector-components of the first principal direction on the right. Blue (resp. red) corresponds to rank 0 (resp. 1) curves.}\label{fig-ecPCA}
\end{figure}

\subsection{Unsupervised learning} 
Unsupervised learning refers to the paradigm in which one seeks to categorize unlabeled data \cite{James2013}. 
Though there are many strategies, we will restrict our exposition to \textsl{dimensionality reduction}. 
The objective is to take data with a high-dimensional representation and to project it down to a lower dimension. 
There are many different techniques to achieve this, ranging from a simple linear  projections to more complicated statistical and topological methods such as $t$-SNE \cite{vanDerMaaten2008} and UMAP \cite{McInnes2018}. 
We will focus our attention on \textsl{principal component analysis} (PCA), which constructs a projection informed by the variance in the dataset.
More precisely, PCA determines a projection onto the space spanned by special vectors along which the data has maximum variance.
These special vectors are called principal components. 

Given an elliptic curve $E$, one may construct a vector of the form $X(E) = ( a_{p_1}(E), a_{p_2}(E), \ldots , a_{p_n}(E))$ where $p_n$ is the $n$th prime. Our high-dimensional data is simply a collection of vectors $X(E)$.
The projection is determined by the variance of $X(E)$ amongst the different $E$, and yields a lower dimensional image that may be calculated by a weighted sum of $a_p(E)$.
The result of two-dimensional PCA is depicted in Figure~\ref{fig-ecPCA}.
The separation in the left frame demonstrates 
a correlation between rank and the projection of the vector $X(E)$ along the direction of maximum variance (the first principal component). 
More surprisingly, in the right frame, we see a structure in the vector-components of this projection. 
This structure is exactly the murmuration phenomenon, but observed through unsupervised machine learning. 
What is remarkable about this approach is that it does not require us to group curves by rank parity, and the murmuration phenomenon emerges simply from the geometry of a high-dimensional point cloud. 
  
\subsection{Supervised learning} 
Supervised learning amounts to a strategy for approximating some function of interest based on a sample of known values \cite{James2013}.
In number theory, a natural function to consider is the one which takes an elliptic curve as its input and produces a rank as its output. 
There is no proven algorithm for calculating the rank function rigorously.
Of course, the BSD conjecture links the rank to an analytic property of the elliptic $L$-function, and that is the starting point for our supervised learning experiments.
In fact, the BSD allows one to estimate the rank of $E$ as a weighted sum of $a_p(E)$.
To further probe this relationship, in \cite{He:2020tkg} we constructed a logistic regression model. Given an elliptic curve $E$, we again consider the vector $X(E)$.
The model predicts 
\[\mathrm{rank}(E) \approx \mathrm{argmax} \ \sigma(X(E) \cdot w + b),\] 
where $\sigma(z)$ is the softmax function generating probability distribution of possible outcomes, and $(w;b)$ are parameters that are learned by minimizing categorical cross-entropy. 
This model was very successful \cite{He:2022pqn}, as indicated in Table~\ref{tab5}.
From an AI perspective, this success indicates that, for curves of different rank, the vectors $X(E)$ are drawn from different distributions. 
To investigate them, a natural first step is to find their centres. This is done by computing the averaging function $m_{\varepsilon,I}(p)$, thus revealing the murmurations from a different perspective.

\begin{table}[h]
\begin{center}
\caption{Distinguishing elliptic curve rank $r_E$ from $a_p$ coefficients using logistic regression.
All data used a random sample of $2.0\times10^4$ curves for each rank, all with conductor in the interval $[1, 1 \times 10^5]$.
Rows 1 - 3 are two-way classifications and row 4 is a three-way classification.
}\label{tab5} 
\begin{tabular}{@{}cccc@{}}
\toprule
$r_E$ & $\lvert$Data$\rvert$= $\#\{E\}$&Precision&Confidence\\
\midrule
$\{0, 1\}$  & 
                     $2.0\times10^4$ ($\times2$) & 0.961 & 0.92 \\

\{0, 2\}& " &  0.996& 0.99\\
 \{1, 2\} &"&  0.999&0.99\\
 \{0, 1, 2\} & $2.0\times10^4$ ($\times3$) & 0.975 & 0.96 \\
\botrule
\end{tabular}
\end{center}
\end{table}

We remark that, had Birch and Swinnerton-Dyer not made their conjecture, an effort that spanned seven years at that time from initial computations, the supervised learning experiment described above would likely have uncovered it naturally.
Indeed, in their first paper on the topic \cite{BIRCHSwinnertonDyer+1963+7+25}, they explain that their motivation is to relate the rank of $E$ to the solution counts $\# E_p$. 
A modern data scientific approach would be to first center the data $\# E_p$, leading naturally to the $a_p(E)$ coefficients as features, before training a logistic regression model to predict the rank from the $a_p(E)$ coefficients. 
Inspecting the learned weights would then reveal that the model is computing a close approximation to the weighted sum of $a_p(E)$ described by BSD conjecture. 

\subsection{Discussion} 
As the experts in attendance at \cite{ICERM} remarked, murmurations could have been spotted in the data available to Birch and Swinnerton-Dyer in 1950s.
In reality, however, murmurations were discovered by an AI-human hybrid.

Mathematicians have used computers to support their research for decades.
Though it has recently entered the public consciousness, machine learning has existed for a similar length of time \cite{5392560}.
What is new is the availability of large databases, which have lowered the barrier to entry for both pure and experimental research in mathematics.
It is obvious that AI is becoming an integral part of mathematical research, in formulating conjectures by pattern recognition, in constructing large language models for mathematical formalism, and in automated proofs and derivations.
From the early experimentations of machine learning in finding short-cuts and new formulae \cite{He:2017aed,He:2018jtw,lample2019deep,udrescu2020ai} to 
the framework in \cite{davies2021advancing} of AI-guidance to conjecture formulation,
it is acknowledged that humans are required for the interpretation and formal proof.

In its current state, it is clear that we are a long way from fully non-human AI mathematicians.
For the pragmatically minded, this raises several questions: What can current AI methodologies say about the deepest questions in mathematics, such as the Riemann Hypothesis or the BSD conjecture?
Can they discover anything new in the patterns for primes, arguably the most elusive of distributions?
Though it was speculated that such discoveries were unlikely \cite{He:2018jtw,kolpakov2023impossibility,wu2023classification}, the current work reports exactly such a thing, the discovery of which was underpinned by AI-assisted experiments and human interpretation of the results. 
In other words, AI-led human intuition and human-led AI exploration should be in tandem for new discoveries in mathematics for the foreseeable future.

\bibliographystyle{unsrt}  
\bibliography{survey-bibliography} 

\end{document}